\documentclass[11pt,reqno]{amsart}
\usepackage{latexsym,mathrsfs}
\usepackage{amssymb,amsfonts,amsmath,mathrsfs,graphicx,color}
\newtheorem{theorem}{Theorem}[section]

\newtheorem{lemma}{Lemma}[section]

\newtheorem{remark}{Remark}[section]

\begin{document}

\title[Gaps between zeros of the Riemann zeta-function]{Gaps between zeros of the Riemann zeta-function}
\author{H. M. Bui and M. B. Milinovich}
\thanks{The paper was commenced while HMB was supported by the Leverhulme Trust RPG-049 at the University of Bristol. MBM is supported in part by NSA Young Investigator grants H98230-15-1-0231 and H98230-16-1-0311.}
\address{School of Mathematics, University of Manchester, Manchester M13 9PL, UK}
\email{hung.bui@manchester.ac.uk}
\address{Department of Mathematics, University of Mississippi, University, MS 38677 USA}
\email{mbmilino@olemiss.edu}

\begin{abstract}
We prove that there exist infinitely many consecutive zeros of the Riemann zeta-function on the critical line whose gaps are greater than $3.18$ times the average spacing. Using a modification of our method, we also show that there are even larger gaps between the multiple zeros of the zeta function on the critical line (if such zeros exist).
\end{abstract}

\allowdisplaybreaks

\maketitle


\section{Introduction}

Let $\zeta(s)$ denote the Riemann zeta-function. In this paper, we investigate the vertical distribution of the nontrivial zeros of $\zeta(s)$, for the most part 
restricting our attention to the zeros located on the critical line $\mathrm{Re}(s)=1/2$. In particular, we show that there are gaps between consecutive zeros on the critical line that are much larger than the average spacing. We also show that there are even larger gaps between the multiple zeros of $\zeta(s)$ (if such zeros exist).

\subsection{Gaps between zeros of the zeta-function}

Let $0<\gamma_1\leq\gamma_2\leq\ldots\leq\gamma_n\leq\ldots$ denote the ordinates of the nontrivial zeros of the Riemann zeta-function in the upper half-plane, and let $t_n$ denote the imaginary part of the $n$th zero of $\zeta(s)$ on the critical line above the real axis. Here, if a zero of the zeta function has multiplicity $m$, then its ordinate is repeated $m$ times in either sequence $\{t_n\}$ or $\{\gamma_n\}$. The Riemann Hypothesis (RH) states that all the nontrivial zeros of the zeta function are on the critical line and so we expect that $\gamma_n=t_n$ for all $n\ge 1$. Our main result applies to the sequence $\{t_n\}$ and is unconditional, but if we assume RH then of course this result applies to the sequence $\{\gamma_n\}$, as well.

It is known that, for $T\geq 10$,
\[
N(T):=\sum_{0<\gamma_n \leq T}1=\frac{T\mathscr{L}}{2\pi}-\frac{T}{2\pi}+O(\mathscr{L}), 
\]
where, here and throughout this paper, we set
\[
\mathscr{L}:=\log\frac{T}{2\pi}.
\]
Therefore the average size of $\displaystyle{\frac{\gamma_{n+1}-\gamma_n}{2\pi/\log \gamma_n}}$ is $1$ as $n\to \infty$, and so defining
\begin{equation*}
\lambda:=\limsup_{n\rightarrow\infty}\frac{\gamma_{n+1}-\gamma_n}{2\pi/\log \gamma_n}\qquad \textrm{and}\qquad\Lambda:=\limsup_{n\rightarrow\infty}\frac{t_{n+1}-t_n}{2\pi/\log t_n},
\end{equation*}
it follows that $\Lambda\geq\lambda\geq1$. Our first result proves that there are gaps between successive ordinates of zeros of $\zeta(s)$ on the critical line that are much larger than the average spacing. 

\begin{theorem}\label{gap1}
We have $\Lambda>3.18$. Consequently, assuming RH, we have $\lambda>3.18$.
\end{theorem}

It is widely believed that the zeros of $\zeta(s)$ are all simple, though this remains an open problem. Modifying our proof of Theorem \ref{gap1}, we prove the existence of even larger gaps between multiple zeros of the zeta-function on the critical line (if such zeros exist). We let $t_n^\star$ denote the imaginary part of the $n$th multiple zero of $\zeta(s)$ on the critical line above the real axis. If there are infinitely many multiple zeros of $\zeta(s)$, we set
\[
\Lambda^\star :=\limsup_{n\rightarrow\infty}\frac{t_{n+1}^\star-t_n^\star}{2\pi/\log t_n^\star}
\]
and otherwise we set $\Lambda^\star =\infty$. Note that the value of $\Lambda^\star$ is unaffected by whether or not we choose to count the sequence $\{t_n^\star\}$ with multiplicity. Trivially, we have $\Lambda^\star \ge 2$ since there are at most $N(T)/2$ multiple zeros in the strip $0<\mathrm{Im}(s)\le T.$ More generally, let $N^s(T)$ denote the number of simple zeros of $\zeta(s)$ on the critical line with imaginary part in the interval $(0,T]$. Then if
\[
N^s(T) \ge \big( C+o(1) \big) N(T)
\]
as $T\to \infty$, it follows that there are at most $(1-C +o(1) ) N(T)/2$ multiple zeros on the critical line up to height $T$ and thus $\Lambda^\star \ge 2/(1-C)$. Since Conrey \cite{Conrey} has shown that $C=2/5$ is admissible, we know that $\Lambda^\star \ge 10/3$. We prove the following stronger lower bound for $\Lambda^\star$.

\begin{theorem}\label{gap2}
 We have $\Lambda^\star>4.05$.
\end{theorem}

We remark that Conrey's result can be slightly improved \cite{BCY}, but this improvement only implies that $\Lambda^\star>3.366$. Theorem \ref{gap2} can be improved conditionally. For instance, the result of Bui and Heath-Brown \cite{BHB} on the proportion of simple zeros of $\zeta(s)$ implies that $\Lambda^\star \ge 27/4$ assuming RH.


\subsection{Previous results and conjectures}

The study of the gaps between the zeros of the Riemann zeta-function is an old problem that has received a great deal of attention. We briefly discuss some of the previous results and conjectures to place Theorem \ref{gap1} in context. 

In 1946, Selberg \cite{S} remarked that he could prove $\Lambda>1$. More recently, Bredberg \cite{B1} proved the quantitative estimate  $\Lambda>2.76$. Bredberg's proof, and our approach in the present paper, are variations of a method of Hall \cite{H} (see also \cite{H3,H2}) who had previously shown that $\Lambda>2.63$. We discuss Hall's method in \textsection \ref{hall}, in particular pointing out the novelties in our approach to proving Theorems \ref{gap1} and \ref{gap2}. We remark that Hall's method has also been adapted to study gaps between zeros of zeta and $L$-functions other than $\zeta(s)$, see \cite{TB2,B1,BHTB,TB1}. 

A different method of Mueller \cite{M} has been used in a number of papers to prove lower bounds for $\lambda$ conditional upon RH and its generalizations, see \cite{B, BMN,CGG3,CGG2, FW,MO,Ng,P}. Our result in Theorem \ref{gap1} that $\lambda>3.18$ assuming RH supersedes all of these previous results. Prior to this paper, the strongest known bounds using Mueller's method was that $\lambda>2.9$ assuming RH and that $\lambda>3.072$ assuming the Generalized Riemann Hypothesis for Dirichlet $L$-functions. These results were established in \cite{B2} and \cite{FW}, respectively. 

It is believed that $\Lambda=\lambda=\infty$. This conjecture is stated by Montgomery \cite{Mo} in his original paper on the pair correlation of the zeros of $\zeta(s)$. Montgomery arrives at this conjecture from the stronger hypothesis that, appropriately normalized, statistics of the nontrivial zeros of the zeta-function should asymptotically behave like the statistics of eigenvalues of large random matrices from the Gaussian Unitary Ensemble (GUE). Indeed this GUE hypothesis suggests that the gaps $\gamma_{n+1}-\gamma_n$ should get as large as $1/\sqrt{\log \gamma_n}$. In this direction, Ben Arous and Bourgade \cite[Section 1.3]{BB} have proposed the more precise conjecture that
\[
\limsup_{n\to \infty} \ (\gamma_{n+1}-\gamma_n) \sqrt{\frac{\log \gamma_n}{32}} = 1.
\]

Unconditionally, Littlewood has shown that $\gamma_{n+1}-\gamma_n = O(1/\log\log\log \gamma_n)$ as $n\to \infty$ while Goldston and Gonek \cite{GG}, sharpening another result of Littlewood, proved that
\[
\limsup_{n\to \infty} \ (\gamma_{n+1}-\gamma_n) \log\log \gamma_n \le \pi
\]
assuming RH.  These results appear to be the best known upper bounds for gaps between consecutive zeros of the zeta-function.


\section{Inequalities, mean value estimates, and numerical calculations}\label{hall}

Modifying an argument of Hall \cite{H}, using classical Wirtinger type inequalities, we reduce the problem of detecting large gaps between zeros of the Riemann zeta-function on the critical line to estimating certain mean-values of $\zeta(s)$ and its derivatives.


\subsection{Wirtinger type inequalities} 

\begin{theorem}\label{W} 
Let $f$ and $f'$ be complex-valued continuous functions on the interval $[a,b]$. 
\begin{enumerate}
\item[(i)] \label{W2}
If $f(a)=f(b)=0$, then
\[
\int_a^b |f(t)|^2  \, \mathrm{d}t \ \le \ \Big(\frac{b-a}{\pi}\Big)^{\!2} \int_a^b |f'(t)|^2 \, \mathrm{d}t.
\]
\item[(ii)] 
If $f(a)=f(b)$ and $\displaystyle{\int_a^b f(t) \, \mathrm{d}t =0}$, then
\[
\int_a^b |f(t)|^2  \, \mathrm{d}t \ \le \ \Big(\frac{b-a}{2\pi}\Big)^{\!2} \int_a^b |f'(t)|^2 \, \mathrm{d}t.
\]
\end{enumerate}
\end{theorem}

Proofs of (i) and (ii) can be found in \cite{HLP}, Theorems 257 and 258, where it is shown that these inequalities hold for functions from $[a,b] \mapsto \mathbb{R}$. The theorem can be extended to complex-valued functions in a straightforward manner by applying the inequalities for real-valued functions to the real and imaginary parts of $f$ separately and then adding. The inequality (i) is sometimes referred to as Wirtinger's inequality in the literature. It is not clear how old these inequalities are or who first proved them. For instance, a proof of (i) was given by Scheeffer \cite{W2} in 1885 and a proof of (ii) was given by Almansi \cite{W1} in 1905.


\subsection{Reduction of Theorem \ref{gap1} to mean value estimates}\label{subsec}
 Suppose, for the sake of contradiction, that 
\begin{equation}\label{assum}
\Lambda\leq \kappa. 
\end{equation}
Let $M(s)$ be a Dirichlet polynomial (chosen to ``amplify" the zeta-function on the critical line) and let
\[
F(t,v,\kappa,M):=e^{ivt \mathscr{L}}\zeta\big(\tfrac{1}{2}+it\big)\zeta\Big(\tfrac{1}{2}+it+i\frac{\kappa\pi}{\mathscr{L}}\Big)M\big(\tfrac{1}{2}+it\big),
\]
where $v\in\mathbb{R}$ is a bounded real number to be chosen later. The factor $e^{ivt \mathscr{L}}$ makes 
$F(t,v,\kappa,M)$ mimic a real-valued function when $T\le t \le 2T$ for a certain choice of $v$ (depending on $M$). In order to simplify a later calculation, we have chosen to use the linear function $\frac{t\mathscr{L}}{2}$ in exponent in place of factor
\[
\theta(t)=\textrm{Im}\bigg(\log\Gamma\Big(\frac{1}{4}+\frac{it}{2}\Big)\bigg)-\frac{(\log\pi)t}{2} 
\]
which appears in the definition of the Hardy $Z$-function, $Z(t)=e^{i\theta(t)} \zeta(\frac{1}{2}+it)$, used in \cite{B1,H}.

Denote the zeros of $F$ in the interval $[T,2T]$ by $\widetilde{t}_1\leq \widetilde{t}_2\leq\ldots\leq \widetilde{t}_N$.  In view of our assumption \eqref{assum}, we have
\[
\widetilde{t}_{n+1}-\widetilde{t}_n\leq \big(1+o(1)\big)\frac{\kappa\pi}{\mathscr{L}}
\]
for $1\le n \le N-1$ as $T\rightarrow\infty$ and so inequality (i) in Theorem \ref{W} implies that
\begin{eqnarray*}
\int_{\widetilde{t}_n}^{\widetilde{t}_{n+1}}|F(t,v,\kappa,M)|^2\,\mathrm{d}t &\leq& \Big( \frac{\widetilde{t}_{n+1}-\widetilde{t}_n}{\pi} \Big)^2 \int_{\widetilde{t}_n}^{\widetilde{t}_{n+1}}|F'(t,v,\kappa,M)|^2\,\mathrm{d}t
\\
&\leq&\big(1+o(1)\big)\frac{\kappa^2}{\mathscr{L}^2}\int_{\widetilde{t}_n}^{\widetilde{t}_{n+1}}|F'(t,v,\kappa,M)|^2\,\mathrm{d}t.
\end{eqnarray*}
Summing over $n$, we derive that
\[
\int_{\widetilde{t}_1}^{\widetilde{t}_{N}} |F(t,v,\kappa,M)|^2 \, \mathrm{d}t\leq\big(1+o(1)\big)\frac{\kappa^2}{\mathscr{L}^2}\int_{\widetilde{t}_1}^{\widetilde{t}_{N}} |F'(t,v,\kappa,M)|^2 \, \mathrm{d}t.
\]
Now, by \eqref{assum}, we see that $\widetilde{t}_1-T$ and $2T-\widetilde{t}_N$ are $\ll 1$. Moreover, our choice of $M(s)$ will ensure that these integrals are $\gg T$ and $|F^{(k)}(t,v,\kappa,M)|^2 \ll_{k,\varepsilon} (|t|+1)^{1-\varepsilon}$ for $\varepsilon>0$, so it follows that
\[
\int_{T}^{2T} |F(t,v,\kappa,M)|^2 \, \mathrm{d}t\leq\big(1+o(1)\big)\frac{\kappa^2}{\mathscr{L}^2} \int_{T}^{2T} |F'(t,v,\kappa,M)|^2 \, \mathrm{d}t.
\]
Therefore, if
\[
h_1(v,\kappa,M):=\limsup_{T\to \infty} \frac{\mathscr{L}^2}{\kappa^2}\frac{\int_{T}^{2T}|F(t,v,\kappa,M)|^2 \, \mathrm{d}t}{\int_{T}^{2T}|F'(t,v,\kappa,M)|^2 \, \mathrm{d}t}>1,
\]
then we have contradicted \eqref{assum} and we may conclude that $\Lambda>\kappa$.


\subsection{Reduction of Theorem \ref{gap2} to mean value estimates}\label{subsec2} We first note that if $a$ and $b$ are multiple zeros of $F$, then
\[
\int_a^b F'(t,v,\kappa,M) \, \mathrm{d}t=0 \quad \text{and} \quad F'(a,v,\kappa,M) = 0 = F'(b,v,\kappa,M).
\]
Therefore, inequality (ii) in Theorem \ref{W} implies that
\[
\int_a^b |F'(t,v,\kappa,M)|^2 \, \mathrm{d}t \le  \left(\frac{b-a}{2\pi}\right)^{\!2} \int_a^b |F''(t,v,\kappa,M)|^2 \, \mathrm{d}t.
\]
Now suppose that $\Lambda^\star \leq \kappa$. Summing over the multiple zeros of $F$ in $[T,2T]$ and arguing as in \textsection \ref{subsec}, it follows that we derive a contradiction if
\[
h_2(v,\kappa,M):=\limsup_{T\to \infty} \frac{4\mathscr{L}^2}{\kappa^2}\frac{\int_{T}^{2T}|F'(t,v,\kappa,M)|^2 \, \mathrm{d}t}{\int_{T}^{2T}|F''(t,v,\kappa,M)|^2 \, \mathrm{d}t}>1,
\]
in which case we can conclude that $\Lambda^\star>\kappa$. Comparing with $h_1(v,\kappa,M)$, we see that  $h_2(v,\kappa,M)$ has an extra factor of 4 in the numerator but the ratio integrals ends up being smaller. Nevertheless, we are able to derive a stronger lower bound for $\Lambda^\star$ than for $\Lambda$.


\subsection{Remarks}

We now point out some of the novelties of our approach. Hall \cite{H} essentially chooses $M(s)=1$ and $v=2$ while  Bredberg \cite{B1} chooses $M(s) =\sum_{h\leq y}1/h^s$ with $y=T^\vartheta$ and $\vartheta<1/11$. We improve upon their results and these choices in number of ways:

\smallskip

\begin{enumerate}
\item[(I)] We choose a more general amplifier of the form
\begin{equation}\label{21}
M(s):=M(s,P) = \sum_{h\leq y}\frac{d_r(h)P[h]}{h^s},
\end{equation}
where $y=T^\vartheta$, $0<\vartheta<1/4$, $r\in\mathbb{N}$, $d_r(h)$ are the coefficients the Dirichlet series of $\zeta(s)^r$, and 
\begin{equation*}
 P[h] := P\Big(\frac{\log y/h}{\log{y}}\Big)
\end{equation*}
for $1 \leq h \leq y$ where $P(x) = \sum_{j\geq0} b_j x^j$ is a certain polynomial. By convention, we set $P[h] = 0$ for $h \geq y$. Note that with this definition we have
\begin{equation}\label{22}
P[h]=\sum_{j\geq 0}\frac{b_j j!}{(\log y)^j}\frac{1}{2\pi i}\int_{(1)}\Big(\frac{y}{h}\Big)^s\frac{\mathrm{d}s}{s^{j+1}}
\end{equation}
for $h \in \mathbb{N}$ (and $y\ne h$ if $j=0$) where here, and throughout the article, the notation $\int_{(c)}$ means $\int_{c-i\infty}^{c+i\infty}$.

\medskip

\item[(II)] In addition to choosing a more general amplifier, we also take the advantage of a longer admissible Dirichlet polynomial in the twisted fourth moment of the zeta function as a consequence of the recent work of Bettin, Bui, Li, and Radziwi\l\l \ \cite{BBLR}. That paper evaluates the integral from $T$ to $2T$ of mean fourth power of the zeta-function on the critical line times the mean square of a Dirichlet polynomial of length $T^\vartheta$ for $\vartheta<1/4$. Bredberg used a result of Hughes and Young \cite{HY} that is valid for $\vartheta<1/11$.

\medskip

\item[(III)] Another novel aspect of our work is that we express our mean value estimates in a more concise and much simpler form. Bredberg's asymptotic formulae took five pages to display while ours are derived from one multiple integral formula (Theorem \ref{smoothed}). Among other things, this helps facilitate numerical calculations.

\medskip

\item[(IV)] Finally, our proof of Theorem \ref{gap2} seems to be the first approach that uses inequality (ii) in Theorem \ref{W} to study the zeros of $\zeta(s)$. 

\medskip

\end{enumerate}


\subsection{A smoothing argument}

To use the result on the twisted fourth moment of the Riemann zeta-function from \cite{BBLR} directly (see Theorem 5.1), we introduce a smooth function $w(t)$ with support in the interval $[1,2]$ and satisfying $w^{(j)}(t)\ll_{j,\varepsilon} T^\varepsilon$ for any $j\geq 0$ and $\varepsilon>0$. For $r\in\mathbb{N}$, we also define the constant $A_r$ via the well known asymptotic formula
\[
\sum_{n\le x} \frac{d_r(n)^2}{n} \sim \frac{A_r   (\log x)^{r^2}}{(r^2)!},
\]
as $x\to\infty$, so
\[
A_r = \prod_{p\text{ prime}} \bigg( \Big(1-\frac{1}{p} \Big)^{r^2} \sum_{\ell=0}^\infty \frac{d_r(p^\ell)^2}{p^\ell}  \bigg).
\]
The smoothed mean values of $|F(t,v,\kappa,M)|^2$, $|F'(t,v,\kappa,M)|^2$, and $|F''(t,v,\kappa,M)|^2$ are given by the following theorem.

\begin{theorem} \label{smoothed}
Suppose that $\vartheta<1/4$. Then, for $T$ large, we have
\begin{equation*}
\int_{-\infty}^{\infty}|F^{(j)}(t,v,\kappa,M)|^2w\Big(\frac tT\Big) \, \mathrm{d}t=\frac{c_j(v,\kappa)A_{r+2}(\log y)^{r^2+4r}\mathscr{L}^{4+2j}}{2(r^2-1)!((r-1)!)^4}\widehat{w}(0)T+O(T\mathscr{L}^{(r+2)^2+2j-1})
\end{equation*}
for $j=0,1,2$, where 
\begin{eqnarray*}
&&\!\!\!\!\!\!\!\!\!\!c_j(v,\kappa)=\mathop{\int}_{\substack{[0,1]^9\\x+x_1+x_2\leq1\\x+x_3+x_4\leq 1}}e^{i\vartheta\kappa\pi\big(x_2-x_4-(x_3-x_4)t_3+(x_1-x_2)t_4\big)-i\kappa\pi\big((1-\vartheta(x_1+x_3))t_1-(1-\vartheta(x_2+x_4))t_2\big)(t_3-t_4)}\\
&&\!\!\!\!\!\!\!\!\!\!\qquad\big(1-\vartheta(x_1+x_3)\big)\big(1-\vartheta(x_2+x_4)\big)\Big(\vartheta(x_1-x_2)+\big(1-\vartheta(x_1+x_3)\big)t_1-\big(1-\vartheta(x_2+x_4)\big)t_2\Big)\\
&&\!\!\!\!\!\!\!\!\!\!\qquad\qquad\Big(\vartheta(x_3-x_4)+\big(1-\vartheta(x_1+x_3)\big)t_1-\big(1-\vartheta(x_2+x_4)\big)t_2\Big)x^{r^2-1}(x_{1}x_2x_3x_4)^{r-1}\nonumber\\
&&\!\!\!\!\!\!\!\!\!\!\qquad\qquad\qquad \Big(v-\vartheta(x+x_1+x_2+x_3+x_4)-\big(1-\vartheta(x_1+x_3)\big)t_1-\big(1-\vartheta(x_2+x_4)\big)t_2\Big)^{2j}\\
&&\!\!\!\!\!\!\!\!\!\!\qquad\qquad\qquad\qquad P(1-x-x_1-x_2)P(1-x-x_3-x_4) \, \mathrm{d}x_1 \, \mathrm{d}x_2 \, \mathrm{d}x_3 \, \mathrm{d}x_4 \, \mathrm{d}x \, \mathrm{d}t_1 \, \mathrm{d}t_2 \, \mathrm{d}t_3 \, \mathrm{d}t_4.
\end{eqnarray*}
\end{theorem}


\subsection{Numerical calculations}

It is a standard exercise to deduce from Theorem \ref{smoothed} the unsmoothed mean-values
\begin{equation*}
\int_{T}^{2T}|F^{(j)}(t,v,\kappa,M)|^2 \, \mathrm{d}t=\frac{c_j(v,\kappa)A_{r+2}(\log y)^{r^2+4r}\mathscr{L}^{4+2j}}{2(r^2-1)!((r-1)!)^4}T+O(T\mathscr{L}^{(r+2)^2+2j-1}).
\end{equation*}
Hence, we deduce from the analysis in \textsection \ref{subsec} and \textsection \ref{subsec2} that
\[
h_1(v,\kappa,M)=\frac{c_0(v,\kappa)}{\kappa^2c_1(v,\kappa)} \qquad\textrm{and}\qquad h_2(v,\kappa,M)=\frac{4c_1(v,\kappa)}{\kappa^2c_2(v,\kappa)}.
\]
A numerical calculation with the values
\[
\vartheta=\tfrac{1}{4},\quad r=1,\quad v=1.26,\quad\textrm{and}\quad P(x)=1-5.8x+6.4x^2
\]
yields
\[
h_1(1.26,3.18,M)>1.0002,
\]
while the values
\[
\vartheta=\tfrac{1}{4},\quad r=1,\quad v=1.25,\quad\textrm{and}\quad P(x)=1-5.2x+5.5x^2
\]
numerically give 
\[
h_2(1.25,4.05,M)>1.0048.
\]
This implies that $\Lambda>3.18$ and $\Lambda^\star > 4.05$ and therefore Theorem \ref{gap1} and Theorem \ref{gap2} follow from Theorem \ref{smoothed}.



\section{A shifted mean value result}

Rather than working directly with the mean squares of $F(t,v,\kappa,M)$, $F'(t,v,\kappa,M)$ and $F''(t,v,\kappa,M)$, we instead consider the shifted mean value
\begin{eqnarray} \label{shifted defn}
I_{\underline{\alpha},\underline{\beta}}(M)&=&\int_{-\infty}^{\infty}\zeta(\tfrac{1}{2}+\alpha_1+it)\zeta(\tfrac{1}{2}+\alpha_2+it)\zeta(\tfrac{1}{2}+\beta_1-it)\zeta(\tfrac{1}{2}+\beta_2-it) \nonumber \\
&&\qquad\qquad\qquad M(\tfrac{1}{2}+\alpha_3+it)M(\tfrac{1}{2}+\beta_3-it)w\Big(\frac tT\Big) \, \mathrm{d}t,
\end{eqnarray}
where the shifts $\alpha_j,\beta_j\ll \mathscr{L}^{-1}$ and the Dirichlet polynomial $M(s)$ is defined in \eqref{21}. Our main goal in the rest of the paper is to prove the following lemma.

\begin{lemma}\label{main}
Suppose $\vartheta<1/4$. Then we have
\[
I_{\underline{\alpha},\underline{\beta}}(M)=\frac{c(\underline{\alpha},\underline{\beta})A_{r+2}(\log y)^{r^2+4r}\mathscr{L}^4}{2(r^2-1)!((r-1)!)^4}\widehat{w}(0)T+O(T\mathscr{L}^{(r+2)^2-1}),
\]
where $c(\underline{\alpha},\underline{\beta})$ is given by
\begin{eqnarray}\label{ca}
&&\!\!\!\!\!\!\!\!\!\!\mathop{\int}_{\substack{[0,1]^9\\x+x_1+x_2\leq1\\x+x_3+x_4\leq 1}}y^{-(\alpha_3+\beta_3)x-\alpha_3( x_1+ x_2)-\beta_3(x_3+x_4)-\beta_1 x_1-\beta_2 x_2-\alpha_1x_3-\alpha_2x_4-(\alpha_2-\alpha_1)(x_3-x_4)t_3-(\beta_2-\beta_1)(x_1-x_2)t_4}\nonumber\\
&&\!\!\!\!\!\!\!\!\!\!\quad(Ty^{-x_1-x_3})^{-(\alpha_1+\beta_1)t_1-(\alpha_2-\alpha_1)t_1t_3-(\beta_2-\beta_1)t_1t_4}(Ty^{-x_2-x_4})^{-(\alpha_2+\beta_2)t_2+(\alpha_2-\alpha_1)t_2t_3+(\beta_2-\beta_1)t_2t_4}\nonumber\\
&&\!\!\!\!\!\!\!\!\!\!\quad\quad\big(1-\vartheta(x_1+x_3)\big)\big(1-\vartheta(x_2+x_4)\big)\Big(\vartheta(x_1-x_2)+\big(1-\vartheta(x_1+x_3)\big)t_1-\big(1-\vartheta(x_2+x_4)\big)t_2\Big)\nonumber\\
&&\!\!\!\!\!\!\!\!\!\!\quad\quad\quad\Big(\vartheta(x_3-x_4)+\big(1-\vartheta(x_1+x_3)\big)t_1-\big(1-\vartheta(x_2+x_4)\big)t_2\Big) x^{r^2-1}(x_{1}x_2x_3x_4)^{r-1}\nonumber\\
&&\!\!\!\!\!\!\!\!\!\!\quad\quad\quad\quad P(1-x-x_1-x_2)P(1-x-x_3-x_4) \, \mathrm{d}x_1 \, \mathrm{d}x_2 \, \mathrm{d}x_3 \, \mathrm{d}x_4 \, \mathrm{d}x \, \mathrm{d}t_1 \, \mathrm{d}t_2 \, \mathrm{d}t_3 \, \mathrm{d}t_4
\end{eqnarray}
uniformly for $\alpha_j,\beta_j\ll\mathscr{L}^{-1}$.
\end{lemma}
We prove this lemma in \textsection \ref{last}. We conclude this section by proving that Theorem 2.3 follows from Lemma \ref{main}. When $j=0$, we have
\[
|F(t,v,\kappa,M)|=\Big|\zeta\big(\tfrac{1}{2}+it\big)\zeta\Big(\tfrac{1}{2}+it+i\frac{\kappa\pi}{\mathscr{L}}\Big)M\big(\tfrac{1}{2}+it\big)\Big|
\]
and hence
\[
c_0(v,\kappa)=c(\underline{\alpha},\underline{\beta})\bigg|_{\substack{\alpha_1=\alpha_3=\beta_1=\beta_3=0\\\alpha_2=i\kappa\pi/\mathscr{L},\beta_2=-i\kappa\pi/\mathscr{L}}}.
\]
In the case $j=1$, we have
\begin{eqnarray*}
&&\!\!\!\!\!\!\!\!\frac{F'(t,v,\kappa,M)}{ie^{2iv\theta(t)}}=v\mathscr{L}\zeta\big(\tfrac{1}{2}+it\big)\zeta\Big(\tfrac{1}{2}+it+i\frac{\kappa\pi}{\mathscr{L}}\Big)M\big(\tfrac{1}{2}+it\big)\\
&&\!\!\!\!\!\!\!\!\qquad\qquad +\Big(\frac{\mathrm{d}}{\mathrm{d}\alpha_1}+\frac{\mathrm{d}}{\mathrm{d}\alpha_2}+\frac{\mathrm{d}}{\mathrm{d}\alpha_3}\Big)\zeta\big(\tfrac{1}{2}+\alpha_1+it\big)\zeta\Big(\tfrac{1}{2}+\alpha_2+it+i\frac{\kappa\pi}{\mathscr{L}}\Big)M\big(\tfrac{1}{2}+\alpha_3+it\big)\bigg|_{\underline{\alpha}=0}\\
&&\!\!\!\!\!\!\!\!\qquad=\mathscr{L}Q\bigg(\frac{1}{\mathscr{L}}\Big(\frac{\mathrm{d}}{\mathrm{d}\alpha_1}+\frac{\mathrm{d}}{\mathrm{d}\alpha_2}+\frac{\mathrm{d}}{\mathrm{d}\alpha_3}\Big)\bigg)\zeta\big(\tfrac{1}{2}+\alpha_1+it\big)\zeta\Big(\tfrac{1}{2}+\alpha_2+it+i\frac{\kappa\pi}{\mathscr{L}}\Big)M\big(\tfrac{1}{2}+\alpha_3+it\big)\bigg|_{\underline{\alpha}=0},
\end{eqnarray*}
where
\[
Q(x)=v+x.
\]
Hence
\begin{eqnarray}\label{abc}
&&\!\!\!\!\!\!\!\!\!\!\!\!\int_{-\infty}^{\infty}|F'(t,v,\kappa,M)|^2w\Big(\frac tT\Big) \, \mathrm{d}t\\
&&=\mathscr{L}^2Q\bigg(\frac{1}{\mathscr{L}}\Big(\frac{\mathrm{d}}{\mathrm{d}\alpha_1}+\frac{\mathrm{d}}{\mathrm{d}\alpha_2}+\frac{\mathrm{d}}{\mathrm{d}\alpha_3}\Big)\bigg)Q\bigg(\frac{1}{\mathscr{L}}\Big(\frac{\mathrm{d}}{\mathrm{d}\beta_1}+\frac{\mathrm{d}}{\mathrm{d}\beta_2}+\frac{\mathrm{d}}{\mathrm{d}\beta_3}\Big)\bigg)I_{\underline{\alpha},\underline{\beta}}(M)\bigg|_{\substack{\alpha_1=\alpha_3=\beta_1=\beta_3=0\\\alpha_2=i\kappa\pi/\mathscr{L},\beta_2=-i\kappa\pi/\mathscr{L}}}.\nonumber
\end{eqnarray}
Similarly, we have
\begin{eqnarray}\label{abcd}
&&\!\!\!\!\!\!\!\!\!\!\!\!\int_{-\infty}^{\infty}|F''(t,v,\kappa,M)|^2w\Big(\frac tT\Big) \, \mathrm{d}t\\
&&=\mathscr{L}^4Q\bigg(\frac{1}{\mathscr{L}}\Big(\frac{\mathrm{d}}{\mathrm{d}\alpha_1}+\frac{\mathrm{d}}{\mathrm{d}\alpha_2}+\frac{\mathrm{d}}{\mathrm{d}\alpha_3}\Big)\bigg)^2Q\bigg(\frac{1}{\mathscr{L}}\Big(\frac{\mathrm{d}}{\mathrm{d}\beta_1}+\frac{\mathrm{d}}{\mathrm{d}\beta_2}+\frac{\mathrm{d}}{\mathrm{d}\beta_3}\Big)\bigg)^2I_{\underline{\alpha},\underline{\beta}}(M)\bigg|_{\substack{\alpha_1=\alpha_3=\beta_1=\beta_3=0\\\alpha_2=i\kappa\pi/\mathscr{L},\beta_2=-i\kappa\pi/\mathscr{L}}}.\nonumber
\end{eqnarray}

We obtain the constants $c_{1}(v,\kappa)$ and $c_{2}(v,\kappa)$ by applying the above differential operators to $c(\underline{\alpha},\underline{\beta})$.  Since $I(\underline{\alpha},\underline{\beta})$ and $c(\underline{\alpha},\underline{\beta})$ are holomorphic with respect to each variable $\alpha_j$ and $\beta_j$ in a small disc centered at $0$, the derivatives appearing in \eqref{abc} and \eqref{abcd} can be expressed using Cauchy's integral formula as integrals of radii $\asymp \mathscr{L}^{-1}$ around the points $\alpha_1=\alpha_3=\beta_1=\beta_3=0,\alpha_2=i\kappa\pi/\mathscr{L},\beta_2=-i\kappa\pi/\mathscr{L}$.  Since the asymptotic formula in Lemma \ref{main} holds uniformly on these contours, each derivative adds a factor that is $O(\mathscr{L})$ to the error term that holds for $I(\underline{\alpha},\underline{\beta})$. Therefore we can use Lemma \ref{main} and \eqref{abc} to prove Theorem \ref{smoothed} in the case $j=1$ with an error of $O(T\mathscr{L}^{(r+2)^2+1})$, and similarly use Lemma \ref{main} and \eqref{abcd} to prove Theorem \ref{smoothed} in the case $j=2$ with an error of $O(T\mathscr{L}^{(r+2)^2+3})$.

To see that applying the above differential operators to $c(\underline{\alpha},\underline{\beta})$ does indeed give the constants $c_{1}(v,\kappa)$ and $c_{2}(v,\kappa)$, note that
\begin{equation*}
\label{eq:Qop}
Q\bigg(\frac{1}{\mathscr{L}}\Big(\frac{d}{d\alpha_1}+\frac{d}{d\alpha_2}+\frac{d}{d\alpha_3}\Big)\bigg)^jX_{1}^{\alpha_1}X_{2}^{\alpha_2}X_{3}^{\alpha_3}=Q\Big(\frac{\log X_1+\log X_2+\log X_3}{\mathscr{L}}\Big)^jX_{1}^{\alpha_1}X_{2}^{\alpha_2}X_{3}^{\alpha_3}.
\end{equation*}
Using this expression and \eqref{ca}, we have
\begin{eqnarray*}
&&\!\!\!\!\!\!\!\!\!\!Q\bigg(\frac{1}{\mathscr{L}}\Big(\frac{d}{d\alpha_1}+\frac{d}{d\alpha_2}+\frac{d}{d\alpha_3}\Big)\bigg)^jQ\bigg(\frac{1}{\mathscr{L}}\Big(\frac{d}{d\beta_1}+\frac{d}{d\beta_2}+\frac{d}{d\beta_3}\Big)\bigg)^jc(\underline{\alpha},\underline{\beta})\\
&&\!\!\!\!\!\!\!\!\!\!\quad=\mathop{\int}_{\substack{[0,1]^9\\x+x_1+x_2\leq1\\x+x_3+x_4\leq 1}}y^{-(\alpha_3+\beta_3)x-\alpha_3( x_1+ x_2)-\beta_3(x_3+x_4)-\beta_1 x_1-\beta_2 x_2-\alpha_1x_3-\alpha_2x_4-(\alpha_2-\alpha_1)(x_3-x_4)t_3-(\beta_2-\beta_1)(x_1-x_2)t_4}\nonumber\\
&&\!\!\!\!\!\!\!\!\!\!\quad\quad(Ty^{-x_1-x_3})^{-(\alpha_1+\beta_1)t_1-(\alpha_2-\alpha_1)t_1t_3-(\beta_2-\beta_1t_1)t_4}(Ty^{-x_2-x_4})^{-(\alpha_2+\beta_2)t_2+(\alpha_2-\alpha_1)t_2t_3+(\beta_2-\beta_1)t_2t_4}\nonumber\\
&&\!\!\!\!\!\!\!\!\!\!\quad\quad\quad\big(1-\vartheta(x_1+x_3)\big)\big(1-\vartheta(x_2+x_4)\big)\Big(\vartheta(x_1-x_2)+\big(1-\vartheta(x_1+x_3)\big)t_1-\big(1-\vartheta(x_2+x_4)\big)t_2\Big)\\
&&\!\!\!\!\!\!\!\!\!\!\quad\quad\quad\quad\Big(\vartheta(x_3-x_4)+\big(1-\vartheta(x_1+x_3)\big)t_1-\big(1-\vartheta(x_2+x_4)\big)t_2\Big) x^{r^2-1}(x_{1}x_2x_3x_4)^{r-1}\nonumber\\
&&\!\!\!\!\!\!\!\!\!\!\quad\quad\quad\quad\quad Q\Big(-\vartheta(x+x_1+x_2+x_3+x_4)-\big(1-\vartheta(x_1+x_3)\big)t_1-\big(1-\vartheta(x_2+x_4)\big)t_2\Big)^{2j}\\
&&\!\!\!\!\!\!\!\!\!\!\quad\quad\quad\quad\quad\quad P(1-x-x_1-x_2)P(1-x-x_3-x_4) \, \mathrm{d}x_1 \, \mathrm{d}x_2 \, \mathrm{d}x_3 \, \mathrm{d}x_4 \, \mathrm{d}x \, \mathrm{d}t_1 \, \mathrm{d}t_2 \, \mathrm{d}t_3 \, \mathrm{d}t_4\nonumber.
\end{eqnarray*}
Theorem \ref{smoothed} now follows by setting $\alpha_1=\alpha_3=\beta_1=\beta_3=0,\alpha_2=i\kappa\pi/\mathscr{L},\beta_2=-i\kappa\pi/\mathscr{L}$ and simplifying.

\section{Two additional lemmas}

\begin{lemma}\label{600}
Let $j \ge 0$, $n\ge 1$ and $r \ge 1$ be integers. Let $y>0$, and let
\begin{equation} \label{firstK}
K_j(\alpha,\beta)=\frac{1}{2\pi i}\int_{(\mathcal{\mathscr{L}}^{-1})}\Big(\frac{y}{n}\Big)^u\zeta^r(1+\alpha+u)\zeta^r(1+\beta+u) \, \frac{\mathrm{d}u}{u^{j+1}}
\end{equation}
with $y\ne n$ if $j=0$. Then we have
\begin{eqnarray*}
K_j(\alpha,\beta)&=&\frac{(\log y/n)^{j+2r}}{((r-1)!)^2j!}\mathop{\int\int}_{\substack{x_1+x_2\leq 1\\0\leq x_1,x_2\leq 1}}\Big(\frac{y}{n}\Big)^{-\alpha x_1-\beta x_2}(x_{1}x_2)^{r-1}(1-x_1-x_2)^j \, \mathrm{d}x_1 \, \mathrm{d}x_2\\
&&\qquad +O\big((\log y)^{j+2r-1}\big)
\end{eqnarray*}
uniformly for $\alpha,\beta\ll (\log y)^{-1}$.
\end{lemma}
\begin{proof}
By a standard application of the residue theorem, we can replace the contour in the integral on the right-hand side of \eqref{firstK} by a small circle with radius $\asymp (\log y)^{-1}$ around the origin plus an error term of size $O(1)$. This integral is trivially bounded by $O\big(\log y)^{j+2r}\big)$. Since
\[
\zeta(1+s) = \frac{1}{s} + O(1) 
\]
for $s$ near zero, taking the first terms in the Laurent series of $\zeta(1+\alpha+u)$ and $\zeta(1+\beta+u)$ gives
\[
K_j(\alpha,\beta)=\frac{1}{2\pi i}\oint q^u\frac{1}{(\alpha+u)^r(\beta+u)^r}\frac{\mathrm{d}u}{u^{j+1}}+O\big((\log y)^{j+2r-1}\big),
\]
where $q=y/n$. We apply the identity
\begin{equation}\label{800}
\frac{1}{(\alpha+u)^r}=\frac{1}{(r-1)!}\int_{1/q}^{1}a^{\alpha+u-1}\Big(\log\frac{1}{a}\Big)^{r-1}\mathrm{d}a+q^{-\alpha-u}\sum_{k=0}^{r-1}\frac{(\log q)^k}{k!(\alpha+u)^{r-k}},
\end{equation}
which is valid for all $\alpha,u\in\mathbb{C}$ and $q>0$, to the above integral, writing it as the sum of $(r+1)$ terms. The last $r$ terms can easily be seen to vanish. Hence
\[
K_j(\alpha,\beta)=\frac{1}{(r-1)!}\int_{1/q}^{1}a^{\alpha-1}\Big(\log\frac{1}{a}\Big)^{r-1}\frac{1}{2 \pi i} \oint (qa)^u\frac{1}{(\beta+u)^r}\frac{\mathrm{d}u}{u^{j+1}}\mathrm{d}a+O\big((\log y)^{j+2r-1}\big).
\] 
We use \eqref{800} again but with the lower boundary of integration at $1/qa$. Again we write the innermost integral as the sum of $(r+1)$ terms where the last $r$ terms vanish. In this way, we derive that
\begin{eqnarray*}
K_j(\alpha,\beta)&=&\frac{1}{((r-1)!)^2}\int_{1/q}^{1}\int_{1/qa}^{1}a^{\alpha-1}b^{\beta-1}\Big(\log\frac{1}{a}\Big)^{r-1}\Big(\log\frac{1}{b}\Big)^{r-1}\frac{1}{2 \pi i} \oint (qab)^u\frac{du}{u^{j+1}} \, \mathrm{d}b \, \mathrm{d}a\\
&&\qquad\qquad+O\big((\log y)^{j+2r-1}\big)\\
&=&\frac{1}{((r-1)!)^2j!}\int_{1/q}^{1}\int_{1/qa}^{1}a^{\alpha-1}b^{\beta-1}\Big(\log\frac{1}{a}\Big)^{r-1}\Big(\log\frac{1}{b}\Big)^{r-1}(\log qab)^j \, \mathrm{d}b \, \mathrm{d}a\\
&&\qquad\qquad+O\big((\log y)^{j+2r-1}\big).
\end{eqnarray*}
Making the variable changes $a \mapsto q^{-x_1}$ and $b \mapsto q^{-x_2}$, we obtain the lemma.
\end{proof}

\begin{lemma}\label{602}
Suppose $f$ is a smooth function. Then we have
\begin{eqnarray*}
\sum_{n\leq y}\frac{d_r(n)}{n^{1+\alpha}}f\Big(\frac{\log y/n}{\log y}\Big)=\frac{(\log y)^r}{(r-1)!}\int_{0}^{1}y^{-\alpha x}x^{r-1}f(1-x) \, \mathrm{d}x+O\big((\log y)^{r-1}\big).
\end{eqnarray*}
\end{lemma}
\begin{proof}
This is a standard exercise in partial summation using the formula
\[
\sum_{n\le y} d_r(n) = \frac{y \, (\log y)^{r-1}}{(r-1)!} + O\big( (\log y)^{r-2} \big),
\]
as $y \to \infty$. See Corollary 4.5 in \cite{BCY}.
\end{proof}

\section{Proof of Lemma \ref{main}}\label{last}

\subsection{Reduction to a contour integral}

We first state the twisted fourth moment of the Riemann zeta-function from \cite{BBLR}. 

\begin{theorem}[Bettin, Bui, Li and Radziwi\l\l]\label{BB}
Let $G(s)$ be an even entire function of rapid decay in any fixed strip $|\emph{Re}(s)|\leq C$ satisfying $G(0)=1$, and let
\begin{equation}\label{Vx}
V(x)=\frac{1}{2\pi i}\int_{(1)}G(s)(2\pi)^{-2s}x^{-s}\frac{\mathrm{d}s}{s}.
\end{equation}
Then, for $T$ large, we have
\begin{align*}
&\sum_{h,k\leq y}\frac{a_h\overline{a_k}}{\sqrt{hk}}\int_{-\infty}^{\infty}\zeta(\tfrac{1}{2}+\alpha_1+it)\zeta(\tfrac{1}{2}+\alpha_2+it)\zeta(\tfrac{1}{2}+\beta_1-it)\zeta(\tfrac{1}{2}+\beta_2-it)\Big(\frac{h}{k}\Big)^{-it}w\Big(\frac{t}{T}\Big) \, \mathrm{d}t\\
&\quad=\sum_{h,k\leq y}\frac{a_h\overline{a_k}}{\sqrt{hk}}\int_{-\infty}^{\infty}w\Big(\frac{t}{T}\Big)\bigg\{Z_{\alpha_1,\alpha_2,\beta_1,\beta_2,h,k}(t)+\Big(\frac{t}{2\pi}\Big)^{-(\alpha_1+\beta_1)}Z_{-\beta_1,\alpha_2,-\alpha_1,\beta_2,h,k}(t)\\
&\quad\quad+\Big(\frac{t}{2\pi}\Big)^{-(\alpha_1+\beta_2)}Z_{-\beta_2,\alpha_2,\beta_1,-\alpha_1,h,k}(t)+\Big(\frac{t}{2\pi}\Big)^{-(\alpha_2+\beta_1)}Z_{\alpha_1,-\beta_1,-\alpha_2,\beta_2,h,k}(t)\\
&\quad\quad\quad+\Big(\frac{t}{2\pi}\Big)^{-(\alpha_2+\beta_2)}Z_{\alpha_1,-\beta_2,\beta_1,-\alpha_2,h,k}(t)+\Big(\frac{t}{2\pi}\Big)^{-(\alpha_1+\alpha_2+\beta_1+\beta_2)}Z_{-\beta_1,-\beta_2,-\alpha_1,-\alpha_2,h,k}(t)\bigg\} \, \mathrm{d}t\\
&\quad\quad\quad\quad+O_\varepsilon(T^{1/2+2\vartheta+\varepsilon}+T^{3/4+\vartheta+\varepsilon})
\end{align*}
uniformly for $\alpha_1,\alpha_2,\beta_1,\beta_2\ll \mathcal{\mathscr{L}}^{-1}$, where $\varepsilon>0$ and
\[
Z_{\alpha,\beta,\gamma,\delta,h,k}(t)=\sum_{hm_1m_2=kn_1n_2}\frac{1}{m_{1}^{1/2+\alpha}m_{2}^{1/2+\beta}n_{1}^{1/2+\gamma}n_{2}^{1/2+\delta}}V\Big(\frac{m_1m_2n_1n_2}{t^2}\Big).
\]
\end{theorem}

\begin{remark}
\emph{To simplify later calculations, it is convenient to prescribe certain conditions on the function $G(s)$. To be precise, we assume that $G(s)$ vanishes at $s=-(\alpha_i+\beta_j)/2$ for $1\leq i,j\leq 2$.}
\end{remark}

Let $\vartheta<1/4$. Recalling the definition of $I_{\underline{\alpha},\underline{\beta}}(M)$ in \eqref{shifted defn}, we write 
\[
I_{\underline{\alpha},\underline{\beta}}(M)=I_1+I_2+I_3+I_4+I_5+I_6+O_\varepsilon(T^{1-\varepsilon})
\]
corresponding to the decomposition in Theorem \ref{BB}. We first estimate $I_1$ using Lemmas \ref{600} and \ref{602}, and then indicate what changes need to be made in our argument to estimate the integrals $I_2,\ldots,I_6.$ Observe that
\[
\begin{split}
I_1&=\sum_{h,k\leq y}\frac{d_r(h)d_r(k)P[h]P[k]}{h^{1/2+\alpha_3}k^{1/2+\beta_3}}\sum_{hm_1m_2=kn_1n_2}\frac{1}{m_{1}^{1/2+\alpha_1}m_{2}^{1/2+\alpha_2}n_{1}^{1/2+\beta_1}n_{2}^{1/2+\beta_2}}\\
&\qquad\qquad\qquad\qquad\int_{-\infty}^{\infty}w\Big(\frac{t}{T}\Big)V\Big(\frac{m_1m_2n_1n_2}{t^2}\Big) \, \mathrm{d}t.
\end{split}
\]
In view of \eqref{22} and \eqref{Vx}, we have
\begin{align} \label{I_1 expression}
I_1=&\sum_{i,j}\frac{b_ib_j i!j!}{(\log y)^{i+j}}\Big(\frac{1}{2\pi i}\Big)^3\int_{-\infty}^{\infty}\int_{(1)}\int_{(1)}\int_{(1)}w\Big(\frac{t}{T}\Big)G(s)\Big(\frac{t}{2\pi}\Big)^{2s}y^{u+v}\\
&\sum_{hm_1m_2=kn_1n_2}\frac{d_r(h)d_r(k)}{h^{1/2+\alpha_3+u}k^{1/2+\beta_3+v}m_{1}^{1/2+\alpha_1+s}m_{2}^{1/2+\alpha_2+s}n_{1}^{1/2+\beta_1+s}n_{2}^{1/2+\beta_2+s}} \frac{\mathrm{d}u}{u^{i+1}}\frac{\mathrm{d}v}{v^{j+1}}\frac{\mathrm{d}s}{s} \, \mathrm{d}t.\nonumber
\end{align}
Since the sum in the integrand is multiplicative, we can express it as an Euler product and then factor out the poles in terms of $\zeta(s)$. In particular, this sum equals
\begin{align}\label{55}
&\sum_{hm_1m_2=kn_1n_2}\frac{d_r(h)d_r(k)}{h^{1/2+\alpha_3+u}k^{1/2+\beta_3+v}m_{1}^{1/2+\alpha_1+s}m_{2}^{1/2+\alpha_2+s}n_{1}^{1/2+\beta_1+s}n_{2}^{1/2+\beta_2+s}}\nonumber\\
&\quad\qquad=A(\underline{\alpha},\underline{\beta},u,v,s)\zeta(1+\alpha_1+\beta_1+2s)\zeta(1+\alpha_1+\beta_2+2s)\zeta(1+\alpha_2+\beta_1+2s)\\
&\quad\quad\qquad\qquad\zeta(1+\alpha_2+\beta_2+2s)\zeta^{r^2}(1+\alpha_3+\beta_3+u+v)\zeta^r(1+\alpha_3+\beta_1+u+s)\nonumber\\
&\quad\quad\quad\qquad\qquad\qquad\zeta^r(1+\alpha_3+\beta_2+u+s)\zeta^r(1+\beta_3+\alpha_1+v+s)\zeta^r(1+\beta_3+\alpha_2+v+s),\nonumber
\end{align}
where $A(\underline{\alpha},\underline{\beta},u,v,s)$ is an arithmetical factor (Euler product) converging absolutely in a product of half-planes containing the origin. 

We first move the $u$ and $v$ contours in \eqref{I_1 expression} to $\textrm{Re}(u)=\textrm{Re}(v)=\delta$, and then move the $s$ contour to $\textrm{Re}(s)=-\delta/2$, where
$\delta> 0$ is some fixed constant such that the arithmetical factor $A(\underline{\alpha},\underline{\beta},u,v,s)$ converges absolutely. In doing so, we cross a pole at $s=0$ and no other singularities of the integrand. Note that the poles at $s=-(\alpha_i+\beta_j)/2$, $1\leq i,j\leq 2$, of the zeta functions are cancelled by the zeros of $G(s)$ and so the integrand is analytic at these points. On the new line of integration we bound the integral by absolute values, giving a contribution of
\[
\ll_\varepsilon T^{1+\varepsilon}y^{2\delta}T^{-\delta}\ll_\varepsilon T^{1-\varepsilon}.
\]
Hence
\begin{eqnarray}\label{40}
I_1&=&\widehat{w}(0)T\zeta(1+\alpha_1+\beta_1)\zeta(1+\alpha_1+\beta_2)\zeta(1+\alpha_2+\beta_1)\zeta(1+\alpha_2+\beta_2)\\
&&\qquad\qquad \sum_{i,j}\frac{b_ib_j i!j!}{(\log y)^{i+j}}J_{i,j}+O_\varepsilon(T^{1-\varepsilon}),\nonumber
\end{eqnarray}
where
\begin{eqnarray*}
J_{i,j}&=&\Big(\frac{1}{2\pi i}\Big)^2\int_{(1)}\int_{(1)}y^{u+v}A(\underline{\alpha},\underline{\beta},u,v,0)\zeta^{r^2}(1+\alpha_3+\beta_3+u+v)\zeta^r(1+\alpha_3+\beta_1+u)\\
&&\qquad\qquad\qquad\zeta^r(1+\alpha_3+\beta_2+u)\zeta^r(1+\beta_3+\alpha_1+v)\zeta^r(1+\beta_3+\alpha_2+v)\frac{\mathrm{d}u}{u^{i+1}}\frac{\mathrm{d}v}{v^{j+1}}.
\end{eqnarray*}

Expressing $\zeta^{r^2}(1+\alpha_3+\beta_3+u+v)$ as an absolutely convergent Dirichlet series and then changing the order of summation and integration, we obtain
\begin{eqnarray}\label{56}
J_{i,j}&=&\sum_{n\leq y}\frac{d_{r^2}(n)}{n^{1+\alpha_3+\beta_3}}\Big(\frac{1}{2\pi i}\Big)^2\int_{(1)}\int_{(1)}\Big(\frac{y}{n}\Big)^{u+v}A(\underline{\alpha},\underline{\beta},u,v,0)\zeta^r(1+\alpha_3+\beta_1+u)\nonumber\\
&&\qquad \zeta^r(1+\alpha_3+\beta_2+u)\zeta^r(1+\beta_3+\alpha_1+v)\zeta^r(1+\beta_3+\alpha_2+v)\frac{\mathrm{d}u}{u^{i+1}}\frac{\mathrm{d}v}{v^{j+1}}.
\end{eqnarray}
Note that here we are able to restrict the sum over $n$ to $n\leq y$ by moving the $u$-integral and the $v$-integral far to the right. We now move the contours of integration to $\textrm{Re}(u)=\textrm{Re}(v)\asymp \mathcal{\mathscr{L}}^{-1}$. Bounding the integrals trivially shows that $J_{i,j}\ll \mathcal{\mathscr{L}}^{i+j+r^2+4r}$. Hence from the Taylor series $A(\underline{\alpha},\underline{\beta},u,v,0)=A(\underline{0},\underline{0},0,0,0)+O(\mathcal{\mathscr{L}}^{-1})+O(|u|+|v|)$, we can replace $A(\underline{\alpha},\underline{\beta},u,v,0)$ by $A(\underline{0},\underline{0},0,0,0)$ in $J_{i,j}$ with an error of size $O(\mathcal{\mathscr{L}}^{i+j+r^2+4r-1})$. We next show that $A(\underline{0},\underline{0},0,0,0)=A_{r+2}$. Letting $\alpha_j=\beta_j=0$ for $j=1,2,3$ and $u=v=s$ in \eqref{55}, we have
\[
\begin{split}
A(\underline{0},\underline{0},s,s,s)&=\zeta(1+2s)^{-(r+2)^2}\sum_{hm_1m_2=kn_1n_2}\frac{d_r(h)d_r(k)}{(hkm_1m_2n_1n_2)^{1/2+s}}\\
&=\zeta(1+2s)^{-(r+2)^2}\sum_{h=k}\frac{d_{r+2}(h)d_{r+2}(k)}{(hk)^{1/2+s}}\\
&=\zeta(1+2s)^{-(r+2)^2}\prod_{p}\sum_{n\geq0}\frac{d_{r+2}(p^n)^2}{p^{n(1+2s)}}.
\end{split}
\]
Hence $A(\underline{0},\underline{0},0,0,0)=A_{r+2}$. The $u$ and $v$ variables in \eqref{56} are now separated so that
\[
J_{i,j}=A_{r+2}\sum_{n\leq y}\frac{d_{r^2}(n)}{n^{1+\alpha_3+\beta_3}}K_{i}(\alpha_3+\beta_1,\alpha_3+\beta_2)K_{j}(\beta_3+\alpha_1,\beta_3+\alpha_2)+O(\mathcal{\mathscr{L}}^{i+j+r^2+4r-1}),
\]
where the function $K_j(\alpha,\beta)$ is defined in Lemma \ref{600}. This lemma implies that
\[
\begin{split}
J_{i,j}&=\frac{A_{r+2}(\log y)^{i+j+4r}}{((r-1)!)^4i!j!}\mathop{\int\int\int\int}_{\substack{x_1+x_2,x_3+x_4\leq 1\\0\leq x_1,x_2,x_3,x_4\leq1}}(x_{1}x_2x_3x_4)^{r-1}(1-x_1-x_2)^i(1-x_3-x_4)^j\\
&\qquad\sum_{n\leq y}\frac{d_{r^2}(n)}{n^{1+\alpha_3+\beta_3}}\Big(\frac{y}{n}\Big)^{-(\alpha_3+\beta_1) x_1-(\alpha_3+\beta_2) x_2-(\beta_3+\alpha_1)x_3-(\beta_3+\alpha_2)x_4}\Big(\frac{\log y/n}{\log y}\Big)^{i+j+4r}\\
&\qquad\qquad\qquad\qquad \mathrm{d}x_1 \, \mathrm{d}x_2 \, \mathrm{d}x_3 \, \mathrm{d}x_4+O(\mathcal{\mathscr{L}}^{i+j+r^2+4r-1}).
\end{split}
\]
Using Lemma \ref{602}, we deduce that
\[
\begin{split}
J_{i,j}&=\frac{A_{r+2}(\log y)^{i+j+r^2+4r}}{(r^2-1)!((r-1)!)^4i!j!}\int_{0}^{1}\mathop{\int\int\int\int}_{\substack{x_1+x_2,x_3+x_4\leq 1\\0\leq x_1,x_2,x_3,x_4\leq1}}\\
&\qquad y^{-(\alpha_3+\beta_3)x-\big((\alpha_3+\beta_1) x_1+(\alpha_3+\beta_2) x_2+(\beta_3+\alpha_1)x_3+(\beta_3+\alpha_2)x_4\big)(1-x)}x^{r^2-1}(1-x)^{i+j+4r}\\
&\qquad\qquad (x_{1}x_2x_3x_4)^{r-1}(1-x_1-x_2)^i(1-x_3-x_4)^j \, \mathrm{d}x_1 \, \mathrm{d}x_2 \, \mathrm{d}x_3 \, \mathrm{d}x_4 \, \mathrm{d}x\\
&\qquad\qquad\qquad+O(\mathcal{\mathscr{L}}^{i+j+r^2+4r-1})\\
&=\frac{A_{r+2}(\log y)^{i+j+r^2+4r}}{(r^2-1)!((r-1)!)^4i!j!}\mathop{\int\int\int\int\int}_{\substack{x+x_1+x_2,x+x_3+x_4\leq 1\\0\leq x,x_1,x_2,x_3,x_4\leq1}}\\
&\qquad y^{-(\alpha_3+\beta_3)x-(\alpha_3+\beta_1) x_1-(\alpha_3+\beta_2) x_2-(\beta_3+\alpha_1)x_3-(\beta_3+\alpha_2)x_4}x^{r^2-1}(x_{1}x_2x_3x_4)^{r-1}\\
&\qquad\qquad (1-x-x_1-x_2)^i(1-x-x_3-x_4)^j \, \mathrm{d}x_1 \, \mathrm{d}x_2 \, \mathrm{d}x_3 \, \mathrm{d}x_4 \, \mathrm{d}x+O(\mathcal{\mathscr{L}}^{i+j+r^2+4r-1}).
\end{split}
\]
Inserting this expression back into \eqref{40}, we conclude that
\begin{eqnarray*}
I_1&=&\frac{A_{r+2}(\log y)^{r^2+4r}\widehat{w}(0)}{(r^2-1)!((r-1)!)^4}\zeta(1+\alpha_1+\beta_1)\zeta(1+\alpha_1+\beta_2)\zeta(1+\alpha_2+\beta_1)\zeta(1+\alpha_2+\beta_2)\\
&&\qquad\mathop{\int\int\int\int\int}_{\substack{x+x_1+x_2,x+x_3+x_4\leq 1\\0\leq x,x_1,x_2,x_3,x_4\leq1}}y^{-(\alpha_3+\beta_3)x-(\alpha_3+\beta_1) x_1-(\alpha_3+\beta_2) x_2-(\beta_3+\alpha_1)x_3-(\beta_3+\alpha_2)x_4}\\
&&\qquad\qquad x^{r^2-1}(x_{1}x_2x_3x_4)^{r-1}P(1-x-x_1-x_2)P(1-x-x_3-x_4) \, \mathrm{d}x_1 \, \mathrm{d}x_2 \, \mathrm{d}x_3 \, \mathrm{d}x_4 \, \mathrm{d}x\\
&&\qquad\qquad\qquad+O(\mathcal{\mathscr{L}}^{r^2+4r-1}).
\end{eqnarray*}

\subsection{Deduction of Lemma \ref{main}} Note that $I_{2}$ is essentially obtained by multiplying $I_{1}$ with $T^{-(\alpha_1+\beta_1)}$ and changing the shifts $\alpha_1\longleftrightarrow-\beta_1$, $I_{3}$ is obtained by multiplying $I_{1}$ with $T^{-(\alpha_1+\beta_2)}$ and changing the shifts $\alpha_1\longleftrightarrow-\beta_2$, $I_{4}$ is obtained by multiplying $I_{1}$ with $T^{-(\alpha_2+\beta_1)}$ and changing the shifts $\alpha_2\longleftrightarrow-\beta_1$, $I_{5}$ is obtained by multiplying $I_{1}$ with $T^{-(\alpha_2+\beta_2)}$ and changing the shifts $\alpha_2\longleftrightarrow-\beta_2$, and $I_{6}$ is obtained by multiplying $I_{1}$ with $T^{-(\alpha_1+\alpha_2+\beta_1+\beta_2)}$ and changing the shifts $\alpha_1\longleftrightarrow-\beta_1$ and $\alpha_2\longleftrightarrow-\beta_2$. Hence
\begin{eqnarray*}
&&I_{\underline{\alpha},\underline{\beta}}(M)=\frac{A_{r+2}(\log y)^{r^2+4r}\widehat{w}(0)}{(r^2-1)!((r-1)!)^4}\mathop{\int\int\int\int\int}_{\substack{x+x_1+x_2,x+x_3+x_4\leq 1\\0\leq x,x_1,x_2,x_3,x_4\leq1}}y^{-(\alpha_3+\beta_3)x-\alpha_3( x_1+ x_2)-\beta_3(x_3+x_4)}\\
&&\qquad\qquad U(\underline{x})x^{r^2-1}(x_{1}x_2x_3x_4)^{r-1}P(1-x-x_1-x_2)P(1-x-x_3-x_4) \, \mathrm{d}x_1 \, \mathrm{d}x_2 \, \mathrm{d}x_3 \, \mathrm{d}x_4 \, \mathrm{d}x\\
&&\qquad\qquad\qquad\qquad+O(\mathcal{\mathscr{L}}^{r^2+4r-1}),
\end{eqnarray*}
where
\[
\begin{split}
U(\underline{x})&=\frac{y^{-\beta_1 x_1-\beta_2 x_2-\alpha_1x_3-\alpha_2x_4}}{(\alpha_1+\beta_1)(\alpha_1+\beta_2)(\alpha_2+\beta_1)(\alpha_2+\beta_2)}-\frac{T^{-(\alpha_1+\beta_1)}y^{\alpha_1 x_1-\beta_2 x_2+\beta_1x_3-\alpha_2x_4}}{(\alpha_1+\beta_1)(-\beta_1+\beta_2)(\alpha_2-\alpha_1)(\alpha_2+\beta_2)}\\
&-\frac{T^{-(\alpha_1+\beta_2)}y^{-\beta_1 x_1+\alpha_1 x_2+\beta_2x_3-\alpha_2x_4}}{(-\beta_2+\beta_1)(\alpha_1+\beta_2)(\alpha_2+\beta_1)(\alpha_2-\alpha_1)}-\frac{T^{-(\alpha_2+\beta_1)}y^{\alpha_2 x_1-\beta_2 x_2-\alpha_1x_3+\beta_1x_4}}{(\alpha_1-\alpha_2)(\alpha_1+\beta_2)(\alpha_2+\beta_1)(-\beta_1+\beta_2)}\\
&-\frac{T^{-(\alpha_2+\beta_2)}y^{-\beta_1 x_1+\alpha_2 x_2-\alpha_1x_3+\beta_2x_4}}{(\alpha_1+\beta_1)(\alpha_1-\alpha_2)(-\beta_2+\beta_1)(\alpha_2+\beta_2)}+\frac{T^{-(\alpha_1+\alpha_2+\beta_1+\beta_2)}y^{\alpha_1 x_1+\alpha_2 x_2+\beta_1x_3+\beta_2x_4}}{(\alpha_1+\beta_1)(\beta_1+\alpha_2)(\beta_2+\alpha_1)(\alpha_2+\beta_2)}.
\end{split}
\]
We write
\[
\begin{split}
\frac{y^{-\beta_1 x_1-\beta_2 x_2-\alpha_1x_3-\alpha_2x_4}}{(\alpha_1+\beta_1)(\alpha_1+\beta_2)(\alpha_2+\beta_1)(\alpha_2+\beta_2)}&=\frac{y^{-\beta_1 x_1-\beta_2 x_2-\alpha_1x_3-\alpha_2x_4}}{(\alpha_1+\beta_1)(-\beta_1+\beta_2)(\alpha_2-\alpha_1)(\alpha_2+\beta_2)}\\
&\qquad-\frac{y^{-\beta_1 x_1-\beta_2 x_2-\alpha_1x_3-\alpha_2x_4}}{(-\beta_1+\beta_2)(\alpha_1+\beta_2)(\alpha_2+\beta_1)(\alpha_2-\alpha_1)}
\end{split}
\]
and
\begin{equation}\label{swap}
\begin{split}
\frac{T^{-(\alpha_1+\alpha_2+\beta_1+\beta_2)}y^{\alpha_1 x_1+\alpha_2 x_2+\beta_1x_3+\beta_2x_4}}{(\alpha_1+\beta_1)(\alpha_1+\beta_2)(\alpha_2+\beta_1)(\alpha_2+\beta_2)}&=\frac{T^{-(\alpha_1+\alpha_2+\beta_1+\beta_2)}y^{\alpha_1 x_1+\alpha_2 x_2+\beta_1x_3+\beta_2x_4}}{(\alpha_1+\beta_1)(-\beta_1+\beta_2)(\alpha_2-\alpha_1)(\alpha_2+\beta_2)}\\
&\qquad-\frac{T^{-(\alpha_1+\alpha_2+\beta_1+\beta_2)}y^{\alpha_1 x_1+\alpha_2 x_2+\beta_1x_3+\beta_2x_4}}{(-\beta_1+\beta_2)(\alpha_1+\beta_2)(\alpha_2+\beta_1)(\alpha_2-\alpha_1)}.
\end{split}
\end{equation}
Notice that we can interchange the roles of $x_1$ with $x_2$, or of $x_3$ with $x_4$, in any term of $U(\underline{x})$ without affecting the value of $I_{\underline{\alpha},\underline{\beta}}(M)$. Applying both changes to the last term on the right-hand side of \eqref{swap}, we can replace $U(\underline{x})$ with the expression
\begin{eqnarray*}
&&\frac{y^{-\beta_1 x_1-\beta_2 x_2-\alpha_1x_3-\alpha_2x_4}}{(-\beta_1+\beta_2)(\alpha_2-\alpha_1)}\bigg(\frac{1-(Ty^{-x_1-x_3})^{-(\alpha_1+\beta_1)}}{\alpha_1+\beta_1}\bigg)\bigg(\frac{1-(Ty^{-x_2-x_4})^{-(\alpha_2+\beta_2)}}{\alpha_2+\beta_2}\bigg)\\
&&\qquad-\frac{y^{-\beta_1 x_1-\beta_2 x_2-\alpha_1x_3-\alpha_2x_4}}{(-\beta_1+\beta_2)(\alpha_2-\alpha_1)}\bigg(\frac{1-(Ty^{-x_2-x_3})^{-(\alpha_1+\beta_2)}}{\alpha_1+\beta_2}\bigg)\bigg(\frac{1-(Ty^{-x_1-x_4})^{-(\alpha_2+\beta_1)}}{\alpha_2+\beta_1}\bigg).
\end{eqnarray*}
Using the integral formula
\begin{equation}\label{int}
\frac{1-z^{-(\alpha+\beta)}}{\alpha+\beta}=(\log z)\int_{0}^{1}z^{-(\alpha+\beta)t} \, \mathrm{d}t,
\end{equation}
we find that
\begin{eqnarray*}
&&\!\!\!\!\!\!\!\!\!\!I_{\underline{\alpha},\underline{\beta}}(M)=\frac{A_{r+2}(\log y)^{r^2+4r}\mathscr{L}^2\widehat{w}(0)}{(r^2-1)!((r-1)!)^4(-\beta_1+\beta_2)(\alpha_2-\alpha_1)}\mathop{\int\int\int\int\int\int\int}_{\substack{x+x_1+x_2,x+x_3+x_4\leq 1\\0\leq x,x_1,x_2,x_3,x_4,t_1,t_2\leq1}}y^{-\beta_1 x_1-\beta_2 x_2-\alpha_1x_3-\alpha_2x_4}\\
&&\!\!\!\!\!\!\!\!\!\!\qquad y^{-(\alpha_3+\beta_3)x-\alpha_3( x_1+ x_2)-\beta_3(x_3+x_4)}(Ty^{-x_1-x_3})^{-(\alpha_1+\beta_1)t_1}(Ty^{-x_2-x_4})^{-(\alpha_2+\beta_2)t_2}\\
&&\!\!\!\!\!\!\!\!\!\!\qquad\qquad \big(1-\vartheta(x_1+x_3)\big)\big(1-\vartheta(x_2+x_4)\big)x^{r^2-1}(x_{1}x_2x_3x_4)^{r-1}\\
&&\!\!\!\!\!\!\!\!\!\!\qquad\qquad\qquad P(1-x-x_1-x_2)P(1-x-x_3-x_4) \, \mathrm{d}x_1 \, \mathrm{d}x_2 \, \mathrm{d}x_3 \, \mathrm{d}x_4 \, \mathrm{d}x \, \mathrm{d}t_1 \, \mathrm{d}t_2\\
&&\qquad-\frac{A_{r+2}(\log y)^{r^2+4r}\mathscr{L}^2}{(r^2-1)!((r-1)!)^4(-\beta_1+\beta_2)(\alpha_2-\alpha_1)}\mathop{\int\int\int\int\int\int\int}_{\substack{x+x_1+x_2,x+x_3+x_4\leq 1\\0\leq x,x_1,x_2,x_3,x_4,t_1,t_2\leq1}}y^{-\beta_1 x_1-\beta_2 x_2-\alpha_1x_3-\alpha_2x_4}\\
&&\!\!\!\!\!\!\!\!\!\!\qquad y^{-(\alpha_3+\beta_3)x-\alpha_3( x_1+ x_2)-\beta_3(x_3+x_4)}(Ty^{-x_2-x_3})^{-(\alpha_1+\beta_2)t_1}(Ty^{-x_1-x_4})^{-(\alpha_2+\beta_1)t_2}\\
&&\!\!\!\!\!\!\!\!\!\!\qquad\qquad \big(1-\vartheta(x_2+x_3)\big)\big(1-\vartheta(x_1+x_4)\big)x^{r^2-1}(x_{1}x_2x_3x_4)^{r-1}\\
&&\!\!\!\!\!\!\!\!\!\!\qquad\qquad\qquad P(1-x-x_1-x_2)P(1-x-x_3-x_4) \, \mathrm{d}x_1 \, \mathrm{d}x_2 \, \mathrm{d}x_3 \, \mathrm{d}x_4 \, \mathrm{d}x \, \mathrm{d}t_1 \, \mathrm{d}t_2+O(\mathcal{\mathscr{L}}^{r^2+4r-1}).
\end{eqnarray*}
Denote the two integrands by $V_1(x,x_1,x_2,x_3,x_4,t_1,t_2)$ and $V_2(x,x_1,x_2,x_3,x_4,t_1,t_2)$, respectively. We note that $I_{\underline{\alpha},\underline{\beta}}(M)$ is unchanged if we swap any of these pairs of variables $x_1\longleftrightarrow x_2$, $x_3\longleftrightarrow x_4$, and $t_1\longleftrightarrow t_2$ in the integrands. Hence we can write
\begin{eqnarray*}
&&\!\!\!\!\!I_{\underline{\alpha},\underline{\beta}}(M)=\frac{A_{r+2}(\log y)^{r^2+4r}\mathscr{L}^2\widehat{w}(0)}{2(r^2-1)!((r-1)!)^4(-\beta_1+\beta_2)(\alpha_2-\alpha_1)}\mathop{\int\int\int\int\int\int\int}_{\substack{x+x_1+x_2,x+x_3+x_4\leq 1\\0\leq x,x_1,x_2,x_3,x_4,t_1,t_2\leq1}}\\
&&\!\!\!\!\!\qquad\qquad\Big(V_1(x,x_1,x_2,x_3,x_4,t_1,t_2)+V_1(x,x_2,x_1,x_4,x_3,t_2,t_1)-V_2(x,x_2,x_1,x_3,x_4,t_1,t_2)\\
&&\!\!\!\!\!\qquad\qquad\qquad\qquad-V_2(x,x_1,x_2,x_4,x_3,t_2,t_1)\Big) \, \mathrm{d}x_1 \, \mathrm{d}x_2 \, \mathrm{d}x_3 \, \mathrm{d}x_4 \, \mathrm{d}x \, \mathrm{d}t_1 \, \mathrm{d}t_2+O(\mathcal{\mathscr{L}}^{r^2+4r-1}).
\end{eqnarray*}
Thus
\begin{eqnarray}\label{If}
&&\!\!\!\!\!\!\!\!\!\!I_{\underline{\alpha},\underline{\beta}}(M)=\frac{A_{r+2}(\log y)^{r^2+4r}\mathscr{L}^2\widehat{w}(0)}{2(r^2-1)!((r-1)!)^4}\mathop{\int\int\int\int\int\int\int}_{\substack{x+x_1+x_2,x+x_3+x_4\leq 1\\0\leq x,x_1,x_2,x_3,x_4,t_1,t_2\leq1}}y^{-(\alpha_3+\beta_3)x-\alpha_3( x_1+ x_2)-\beta_3(x_3+x_4)}\nonumber\\
&&\!\!\!\!\!\!\!\!\!\!\qquad W(\underline{x},t_1,t_2)\big(1-\vartheta(x_1+x_3)\big)\big(1-\vartheta(x_2+x_4)\big)x^{r^2-1}(x_{1}x_2x_3x_4)^{r-1}\\
&&\!\!\!\!\!\!\!\!\!\!\qquad\qquad P(1-x-x_1-x_2)P(1-x-x_3-x_4) \, \mathrm{d}x_1 \, \mathrm{d}x_2 \, \mathrm{d}x_3 \, \mathrm{d}x_4 \, \mathrm{d}x \, \mathrm{d}t_1 \, \mathrm{d}t_2+O(\mathcal{\mathscr{L}}^{r^2+4r-1}),\nonumber
\end{eqnarray}
where
\[
\begin{split}
W(\underline{x},t_1,t_2)&=\frac{1}{(-\beta_1+\beta_2)(\alpha_2-\alpha_1)}\bigg(y^{-\beta_1 x_1-\beta_2 x_2-\alpha_1x_3-\alpha_2x_4}(Ty^{-x_1-x_3})^{-(\alpha_1+\beta_1)t_1}(Ty^{-x_2-x_4})^{-(\alpha_2+\beta_2)t_2}\\
&\qquad +y^{-\beta_1 x_2-\beta_2 x_1-\alpha_1x_4-\alpha_2x_3}(Ty^{-x_2-x_4})^{-(\alpha_1+\beta_1)t_2}(Ty^{-x_1-x_3})^{-(\alpha_2+\beta_2)t_1}\\
&\qquad\qquad-y^{-\beta_1 x_2-\beta_2 x_1-\alpha_1x_3-\alpha_2x_4}(Ty^{-x_1-x_3})^{-(\alpha_1+\beta_2)t_1}(Ty^{-x_2-x_4})^{-(\alpha_2+\beta_1)t_2}\\
&\qquad\qquad\qquad-y^{-\beta_1 x_1-\beta_2 x_2-\alpha_1x_4-\alpha_2x_3}(Ty^{-x_2-x_4})^{-(\alpha_1+\beta_2)t_2}(Ty^{-x_1-x_3})^{-(\alpha_2+\beta_1)t_1}\bigg).
\end{split}
\]
Using \eqref{int} again we see that
\[
\begin{split}
W(\underline{x},t_1,t_2)&=y^{-\beta_1 x_1-\beta_2 x_2-\alpha_1x_3-\alpha_2x_4}(Ty^{-x_1-x_3})^{-(\alpha_1+\beta_1)t_1}(Ty^{-x_2-x_4})^{-(\alpha_2+\beta_2)t_2}\\
&\qquad\quad\bigg(\frac{1-\big(y^{x_3-x_4}(Ty^{-x_1-x_3})^{t_1}(Ty^{-x_2-x_4})^{-t_2}\big)^{-(\alpha_2-\alpha_1)}}{\alpha_2-\alpha_1}\bigg)\\
&\qquad\qquad\quad\quad\bigg(\frac{1-\big(y^{x_1-x_2}(Ty^{-x_1-x_3})^{t_1}(Ty^{-x_2-x_4})^{-t_2}\big)^{-(\beta_2-\beta_1)}}{\beta_2-\beta_1}\bigg)\\
&=\mathscr{L}^2\Big(\vartheta(x_1-x_2)+t_1\big(1-\vartheta(x_1+x_3)\big)-t_2\big(1-\vartheta(x_2+x_4)\big)\Big)\\
&\qquad\quad\Big(\vartheta(x_3-x_4)+t_1\big(1-\vartheta(x_1+x_3)\big)-t_2\big(1-\vartheta(x_2+x_4)\big)\Big)\\
&\qquad\qquad\quad\quad y^{-\beta_1 x_1-\beta_2 x_2-\alpha_1x_3-\alpha_2x_4}(Ty^{-x_1-x_3})^{-(\alpha_1+\beta_1)t_1}(Ty^{-x_2-x_4})^{-(\alpha_2+\beta_2)t_2}\\
&\qquad\qquad\qquad\quad\quad\quad\int_{0}^{1}\int_{0}^{1}\big(y^{x_3-x_4}(Ty^{-x_1-x_3})^{t_1}(Ty^{-x_2-x_4})^{-t_2}\big)^{-(\alpha_2-\alpha_1)t_3}\\
&\qquad\qquad\qquad\qquad\quad\quad\quad\quad\big(y^{x_1-x_2}(Ty^{-x_1-x_3})^{t_1}(Ty^{-x_2-x_4})^{-t_2}\big)^{-(\beta_2-\beta_1)t_4} \, \mathrm{d}t_3 \, \mathrm{d}t_4.
\end{split}
\]
Using this expression in \eqref{If}, we obtain the lemma.

\end{document}